\documentclass[11pt]{article}

\usepackage{graphicx}
\usepackage{amsmath,amsfonts,amssymb}
\usepackage{citesort}
\usepackage{url}
\usepackage{amsmath}
\usepackage{graphicx}
\usepackage{multirow}
\usepackage{citesort}
\usepackage{epsfig}
\usepackage{color}
\def\tto{\;{\lower 1pt \hbox{$\rightarrow$}}\kern -10pt
\hbox{\raise 2pt \hbox{$\rightarrow$}}\;}

\def\Tilde{\widetilde}
\def\Bar{\overline}
\def\ra{\rangle}
\def\la{\langle}
\def\ve{\varepsilon}
\def\B{I\!\!B}
\def\h{\hfill\square}
\def\R{I\!\!R}
\def\N{I\!\!N}
\def\ox{\bar{x}}

\def\co{\mbox{\rm co}\,}

\def\epi{\mbox{\rm epi}\,}

\def\h{\hfill\square}

\def\ph{\varphi}
\def\emp{\emptyset}

\def \N{I\!\!N}

\def\Finf{F_{\infty}}
\newcounter{lk}

\begin{document}
\begin{center}
\vspace*{0.3in} {\bf MINIMAL TIME FUNCTIONS AND THE SMALLEST INTERSECTING BALL PROBLEM GENERATED BY UNBOUNDED DYNAMICS}\\[2ex]
Nguyen Mau Nam\footnote{Department of Mathematics, University of
Texas-Pan American, TX 78539, USA, email: nguyenmn@utpa.edu. The
research of this author was partially supported by the Simons
Foundation under grant \#208785.}, Cristina
Villalobos\footnote{Department of Mathematics, The University of
Texas--Pan American, Edinburg, TX 78539--2999, USA (email:
mcvilla@utpa.edu).}, and Nguyen Thai An \footnote{Department of
Mathematics, College of Education, Hue University, 32 Leloi Hue,
Vietnam (email: thaian2784@gmail.com).}
\end{center}

\small{\bf Abstract:}  The \emph{smallest enclosing circle problem}
introduced in the 19th century by J. J. Sylvester \cite{syl} asks for
the circle of smallest radius enclosing a given set of finite points
in the plane. An extension of the smallest enclosing circle problem
called the \emph{smallest intersecting ball problem} was considered
in \cite{mnv,naj}: given a finite number of nonempty closed subsets
of a normed space, find a ball with the smallest radius that
intersects all of the sets. In this paper we initiate the study of
minimal time functions generated by unbounded dynamics and discuss
their applications to extensions of the smallest intersecting ball
problem. This approach continues our effort in applying convex and
nonsmooth analysis to the well-established field of facility
location.

\medskip
\vspace*{0,05in} \noindent {\bf Key words.} Minimal time functions,
subdifferential, subgradient method, the smallest intersecting ball
problem.


\newtheorem{theorem}{Theorem}[section]
\newtheorem{proposition}{Proposition}[section]
\newtheorem{lemma}{Lemma}[section]
\newtheorem{corollary}{Corollary}[section]
\newtheorem{definition}{Definition}[section]
\newtheorem{example}{Example}[section]
\newtheorem{remark}{Remark}[section]

\normalsize

\section{Introduction and Preliminaries}
\setcounter{equation}{0} Given a finite number of points on the
plane, the smallest enclosing circle problem asks for the smallest
circle that encloses all of the points. This celebrated problem was
introduced in the 19th century by the English mathematician J. J.
Sylvester \cite{syl} and is a simple example of continuous facility
location problems. The reader is referred to \cite{ljc,frank,wel}
and the references therein for more recent results involving this
problem in both theoretical and numerical aspects.\vspace*{0.05in}

Let $X$ be a normed space and let $F \subseteq X $ be a nonempty
closed convex subset of $X$. Given nonempty closed target sets $
\Omega_i$ for $i=1,\ldots,m$ and a nonempty closed constraint set
$\Omega \subseteq X$, our goal in this paper is to study the
following extended version of the smallest intersecting ball
problem: find a point $\ox\in\Omega$ and the smallest radius $r\geq
0$ (if they exist) such that the set $x+rF$ intersects all of the
target sets. It is obvious that when $F$ is the Euclidean closed
unit ball of the plane, the target sets under consideration are
singletons, and the constraint set is the whole plane, this problem
reduces to the smallest enclosing circle problem.\vspace*{0.05in}

The so-called minimal time function below plays a crucial role in the study of the smallest intersecting ball problem. For a nonempty subset $Q \subseteq X$, define
\begin{equation}\label{MT}
\mathcal{T}^F_Q(x)=\inf\{t\geq 0: (x+tF)\cap Q\neq \emptyset\},
\end{equation}
which signifies the minimal time for the point $x$ to reach the
target set $Q$ following the dynamic $F$. General and generalized
differentiation properties of the minimal time function (\ref{MT})
have been studied extensively in the literature; see, e.g.
\cite{cowo,cowo1,heng,bmn10,mn10} and the references therein. When $F$ is
the closed unit ball of $X$, the minimal time function becomes the
familiar distance function
\begin{equation*}
d(x; Q)=\inf\{ ||x-q||: q\in Q\}.
\end{equation*}

In \cite{cowo,cowo1,heng,cgm}, the minimal time function (\ref{MT})
in which $F$ is a closed bounded convex set that contains the origin
as an interior point is considered. Further extensions to the case
where the origin is not necessarily an interior point of $F$ have
been considered in \cite{bmn10,jh}. However, to the best of our
knowledge, the minimal time function (\ref{MT}) with $F$ being
unbounded has not been addressed in the literature. In this paper,
we are going to initiate the study of the minimal time function
(\ref{MT}) for the case where $F$ is unbounded and consider applications to the corresponding extended version of the smallest
intersecting ball problem. \vspace*{0.05in}

Throughout the paper, we  make the following standing assumption, unless otherwise specified:
\begin{center}
\emph{$X$ is a normed space and $F \subseteq X$ is a nonempty closed convex set.}
\end{center}

Under natural assumptions, the smallest intersecting ball problem
generated by the dynamic $F$ with the target sets $\Omega_i$ for
$i=1,\ldots,m$ and the constraint set $\Omega$ can be modeled as the
following optimization problem:

\begin{equation}\label{SIB}
\mbox{minimize } \mathcal{T}(x) \mbox{ subject to } x\in \Omega,
\end{equation}
where
\begin{equation*}
\mathcal{T}(x)=\max\{ \mathcal{T}^F_{\Omega_i}(x): i=1,\ldots,m\}.
\end{equation*}

The functions involved in this problem are nondifferentiable with extended-real values in general, and they are convex when all of the sets therein are convex. Thus it is natural to use nonsmooth analysis and convex analysis in particular to study the problem.\vspace*{0.05in}

Given a convex function $\ph: X\to (-\infty, \infty]$ and a point $\ox\in \mbox{dom }\ph$, the subdifferential in the sense of convex analysis of the function $\ph$ at $\ox$ is defined by
\begin{equation*}
\partial \ph(\ox)=\{x^*\in X^*: \la x^*, x-\ox\ra \leq \ph(x)-\ph(\ox) \mbox{ for all }x\in X\},
\end{equation*}
where the set $X^*$ is the dual of $X$ and the domain of the function $\ph$ is
\begin{equation*}
\mbox{dom }\ph=\{x\in X: \ph(x)<\infty\}.
\end{equation*}
An element $x^*\in \partial \ph(\ox)$ is called a subgradient of the function $\ph$ at $\ox$. It is well known that subdifferentials of convex functions play an important role in convex optimization in both theoretical and numerical aspects. For instance, the generalized Fermat rule
\begin{equation*}
\ph \mbox{ has an absolute minimum at }\ox \mbox{ if and only if }0\in \partial \ph(\ox)
\end{equation*}
allows us to use subdifferentials as a tool to solve nonsmooth convex optimization problems. Under natural assumptions, the iterative sequence
\begin{equation*}
x_{k+1}=x_k-\alpha_k x^*_k, \mbox{ with a chosen }x_1\in X \mbox{ and }x^*_k\in \partial \ph(x_k),
\end{equation*}
converges to a minimizer of the function $\ph$; see, e.g., \cite{bert} and the references therein.\vspace*{0.05in}

Due to the involvement of the minimal time function (\ref{MT}) to the smallest intersecting ball problem, knowing subdifferential properties of this function provides important tools for the study of the problem in the convex setting. In fact, this has been shown by our recent papers \cite{naj,mnft,mnj2} that generalized differentiation of the minimal time function (\ref{MT}) plays an important role in the study of \emph{set facility location problems} including the smallest intersecting ball problem. It turns out that subdifferential properties of the minimal time function (\ref{MT}) have close connections to the \emph{asymptotic cone} of $F$ and the normal cone to $Q$ at $\ox\in Q$ defined by
\begin{equation*}
N(\ox; Q)=\{x^*\in X^*: \la x^*, x-\ox\ra \leq 0\mbox{ for all }x\in  Q\}.
\end{equation*}

In this paper, we review important properties of asymptotic cones for convex sets in Section 2. Section 3 is devoted to general properties of the minimal time function (\ref{MT}), while generalized differentiation of this function in the convex case is addressed in Section 4. Finally, in Section 5, we provide applications of the results obtained to numerically solving the optimization problem (\ref{SIB}) and its connection to the smallest intersecting ball problem. We review the literature and provide self-contained proofs for most of the results for the convenience of the reader.\vspace*{0.05in}

Throughout the paper, $X^*$ denotes the dual space of $X$ and $\la x^*,x\ra=x^*(x)$ is the dual pair between an element $x^*\in X^*$ and an element $x\in X$.

\section{Asymptotic Cones}
\setcounter{equation}{0}

In this section we review some facts about asymptotic cones that will be used in our results. For $x \in F$, the asymptotic cone of $F$ at $x$ is defined by
\begin{align*}
F_\infty(x) &= \{ d \in X : x + td \in F \mbox{ for all } t >0 \}\\
&=\{d\in X: x+td\in F \mbox{ for all }t\geq 0\}.
\end{align*}
Another equivalent definition of $F_\infty(x)$ is given as
$$ F_\infty(x)= \bigcap_{t>0} \frac{F-x}{t}, $$
which shows that $F_\infty(x)$ is a cone and contains the origin, 0.  In addition, $F_\infty(x)$ is closed and convex since $F$ is a closed convex set and
the intersection of closed (convex) sets is closed (convex). Thus $F_\infty(x)$ is a closed convex cone.

\begin{proposition} \label{p1} For all $x_1, x_2 \in F$, we have
\begin{equation*}
F_{\infty}(x_1) = F_{\infty}(x_2),
\end{equation*}
that is, the asymptotic cone does not depend on $x \in F$.  Thus we will denote the asymptotic cone as $F_{\infty}$.
\end{proposition}
\noindent{\bf Proof. }For the proof of this statement in finite dimensions,
the reader is referred \cite{hl}. The proof remains valid in
infinite dimensions, and we include the details here for the
convenience of the reader. Fix $x_1$ and $x_2$ in $F$. It suffices
to show that
\begin{equation*}
F_\infty(x_1) \subseteq F_\infty(x_2).
\end{equation*}
Let $d\in F_\infty(x_1)$. Fix any $t>0$, we will show that $x_2+td\in F$. Consider the sequence
\begin{equation*}
x_n=x_1+td +\left(1-\dfrac{1}{n}\right)(x_2-x_1),
\end{equation*}
which can be written as
\begin{equation*}
x_n=\dfrac{1}{n}(x_1 +ntd) + \left(1-\dfrac{1}{n}\right)x_2.
\end{equation*}
We see that $x_n\in F$ for every $n$ because $d \in F_\infty(x_1)$ and $F$ is convex. In addition, $x_n\to x_2+td$. Thus $x_2+td\in F$ since $F$ is closed. Therefore, $d\in F_\infty(x_2)$. $\h$\vspace*{0.05in}

The following corollary is a direct consequence of Proposition \ref{p1}.
\begin{corollary} \label{c1} Suppose $F$ contains the origin. Then
\begin{equation*}
 F_{\infty} = \bigcap_{t >0} tF.
\end{equation*}
\end{corollary}
\noindent{\bf Proof. }Through the use of definitions of $F_{\infty}$ and Proposition \ref{p1}, we have
\begin{equation*}
F_{\infty} = F_{\infty}(0)= \bigcap_{\hat t >0} \frac{F-0}{\hat t} = \bigcap_{\hat t >0} \frac{1}{\hat t} F
= \bigcap_{t>0} tF
\end{equation*}
 where $t = 1/\hat t$ for $\hat t>0$.$\h$

\begin{proposition}\label{c2}
The following are equivalent:\\[1ex]
{\rm (i)} $d \in F_{\infty}$.\\
{\rm (ii)} There exists a sequence $\{ t_n\} \subseteq [0,\infty)$  such that $t_n \rightarrow 0$ and
 a sequence $\{f_n\}\subseteq F$ with $t_n f_n \rightarrow d$.\\
{\rm (iii)} There exists a sequence $\{ t_n\}\subseteq [0,\infty)$  such that $t_n \rightarrow 0$ and
 a sequence $\{f_n\} \subseteq F$ with $t_n f_n \xrightarrow{w} d$.
\end{proposition}
\noindent{\bf Proof. }Let us prove ``(i)$\Rightarrow$ (ii)". Suppose $d \in \Finf$ and fix $\bar x\in F$. Then
\begin{equation*}
\bar x+td\in F \mbox{ for all } t>0.
\end{equation*}
In particular,
\begin{equation*}
\bar x +nd\in F \mbox{ for all }n\in \N.
\end{equation*}
For each $n$, there exists $f_n\in F$ such that
\begin{equation*}
\bar x+ nd=f_n \mbox{ or equivalently } \dfrac{1}{n}\bar x+d=\dfrac{1}{n} f_n.
\end{equation*}
Let $t_n=\dfrac{1}{n}$. We easily see that $t_nf_n \rightarrow  d$ as $n\to\infty$.

The implication ``(ii)$\Rightarrow$ (iii)" is obvious since the strong convergence in $X$ implies the weak convergence. We finally prove ``(iii)$\Rightarrow$ (i)". Assume that there exist sequences $\{t_n\} \subseteq [0,\infty)$  and $\{f_n\} \subseteq F $  such that $t_n \rightarrow 0$ and $t_n f_n \xrightarrow{w} d$.  Fix any $x\in F$, and let $t>0$.  We will prove that $d \in \Finf$ by showing that
\begin{equation*}
x+td\in F\mbox{ for all }t>0.
\end{equation*}
One has $0\leq t\cdot t_n<1$ when $n$ is sufficiently large and $t \cdot t_n\to 0$. Thus
\begin{equation*}
 (1-t \cdot t_n)x+ t \cdot t_nf_n \xrightarrow{w} x+td \mbox{ as } n\to \infty.
\end{equation*}
 The elements $(1-t \cdot t_n)x+ t \cdot t_nf_n$ are in $F$ because $F$ is convex, and  since $F$ is weakly closed, $x+td\in F$. Therefore, $d\in F_\infty$. $\h$

\begin{definition} The set $F$ is called sequentially compact at $0$ if for any sequence $\{f_n\} \subseteq F$ such that $f_n\xrightarrow{w}0$, one has $f_n\to 0$.
\end{definition}
This property holds obviously in finite dimensions. Moreover, if $X$ is a Hilbert space and $F$ is a nonempty closed convex subset of the \emph{Fr\'echet normal cone} to a set that is \emph{sequentially normally compact} (see \cite{mor06a}) at the reference point, then $F$ is automatically sequentially compact at $0$. In infinite dimensions, it is not applicable to sets with nonempty interior.

\begin{proposition}
If $F$ is bounded, then $F_\infty = \{ 0\}$. The converse holds if we suppose further that $X$ is a reflexive Banach space and $F$ is sequentially compact at the origin with $0\in F$ or $X$ is finite dimensional.
\end{proposition}
\noindent{\bf Proof. }Suppose $F$ is bounded and $d\in F_\infty$. By Proposition \ref{c2}, there exist  sequences $\{t_n\}$ such that $t_n\to 0^+$ and $\{f_n\}\subseteq F$ such that
\begin{equation*}
t_nf_n\rightarrow d \mbox{ as } n\to \infty.
\end{equation*}
Since $F$ is bounded, the sequence $\{t_nf_n\}$ converges to $0$, so $d=0$.

To prove the converse, we suppose by contradiction that $F$ is unbounded and $F_{\infty}=\{0\}$, while $F$ is sequentially compact at the origin with $0\in F$ or $X$ is finite dimensional. Then there exits a sequence $\{x_n\}\subseteq F$ with $||x_n||\to\infty$. The sequence defined by
\begin{equation*}
d_n=\dfrac{x_n}{||x_n||}
\end{equation*}
is bounded. So it has a subsequence (without relabeling) such that $d_n\xrightarrow{w}d \in F$. Suppose first that $F$ is sequentially compact at the origin with $0\in F$. Since $0\in F$ and $F$ is convex, $d_n\in F$ when $n$ is sufficiently large. The sequentially compactness of $F$ at $0$ implies $d\neq 0$. Indeed, if $d=0$, then $d_n\to 0$, but this is a contradiction since $||d_n||=1$ for every $n$. Fix any $x\in F$, and let $t>0$. For sufficiently large $k$, one has
\begin{equation*}
u_k=\left(1-\dfrac{t}{||x_k||}\right)x+ \dfrac{t}{||x_k||}x_k\in F
\end{equation*}
since $F$ is convex.  Since $u_k\xrightarrow{w} x+td$ and $F$ is weakly closed, $x+td\in F$. Thus $d\in F_\infty$, which is a contradiction. The proof is similar in the case $X$ is finite dimensional. $\h$

\section{General Properties of Minimal Time Functions with Unbounded Dynamics}
\label{s:SCNC}

In this section we study general properties of the minimal time function (\ref{MT}) that will be used in the next sections. These properties are also of independent interest.

It is obvious that if $x\in Q$, then $\mathcal{T}^F_Q(x)=0$. Moreover, $\mathcal{T}^F_Q(x)$ is finite if there exists $t\geq 0$ such that
\begin{equation*}
(x+tF)\cap Q\neq \emptyset.
\end{equation*}
This holds in particular when $0\in \mbox{int }F$. If no such $t$ exists, then $\mathcal{T}^F_Q(x)=\infty$. Thus the minimal time function (\ref{MT}) is an extended real-valued function in general.

\begin{theorem}\label{t1} Suppose one of the following:\\
{\rm (i)} $X$ is a reflexive Banach space,  $Q$ is a bounded weakly closed set, and $F$ contains the origin $0$. \\
{\rm (ii)} $Q$ is a compact set and $F$ contains the origin $0$.\\[1ex]
Then
$$\mathcal{T}^F_Q(x) = 0  \mbox{ if and only if } x \in Q - \Finf. $$
\end{theorem}
\noindent{\bf Proof. }Let us prove the theorem under assumption (i). Suppose $\mathcal{T}^F_Q(x)=0$. Then for every $n\in \N$, there exists $t_n\geq 0$ such that
\begin{equation*}
(x+t_nF)\cap  Q\neq \emptyset,
\end{equation*}
and $\{t_n\}$ converges to $0$. Then there exist $\omega_n\in  Q$ and $f_n\in F$ for every $n \in \N$ with
\begin{equation*}
x+t_nf_n=\omega_n.
\end{equation*}
Since $Q$ is bounded and weakly closed in a reflexive Banach space, there exists a subsequence (without relabeling) such that the sequence $\{\omega_n\}$ converges weakly to $\omega\in Q$. Thus
\begin{equation*}
t_nf_n\xrightarrow{w} \omega-x \mbox{ as } n\to \infty.
\end{equation*}
It follows from Proposition \ref{c2} that $\omega-x\in F_\infty$, and hence $x\in Q-F_\infty$.

Let us prove the opposite implication. Fix any $x\in Q-F_\infty$. Then there exist $\omega\in Q$ and $d\in F_\infty$ such that
\begin{equation*}
x=\omega-d.
\end{equation*}
Since $0\in F$ and $d \in F_{\infty}$, then by Corollary \ref{c1}, $n(\omega-x)=nd\in F$ for all $n\in \N$. Thus
\begin{equation*}
\omega-x\in \dfrac{1}{n} F\mbox{ for all } n\in \N.
\end{equation*}
This implies $(x+\dfrac{1}{n} F ) \cap Q\neq \emptyset$ for all $n\in \N$. Thus $0\leq \mathcal{T}^F_Q(x)\leq \dfrac{1}{n}$ for all $n\in \N$. Therefore, $\mathcal{T}^F_Q(x)= 0$. The proof of the theorem under  assumption (ii) is similar. $\h$

\begin{remark}\label{r1}{\rm The proof of Theorem \ref{t1} shows that, if $Q$ is a nonempty set (not necessarily bounded) and $0\in F$, then $\mathcal{T}^F_Q(x) =  0 $ whenever $x\in Q-F_\infty$.}
\end{remark}

The following example is a demonstration of Theorem \ref{t1}.

\begin{example}{\rm Let $F=\R\times [-1,1]\subset \R^2$ and let $ Q$ be the disk with center at $(1,0)$ and radius $r=1$. Then using the definition of $F_\infty$ in Corollary \ref{c1} we obtain
\begin{equation*}
F_\infty=\R\times \{0\},
\end{equation*}
and $ Q-F_\infty=\R\times [-1,1]$. Thus $\mathcal{T}^F_ Q(x)=0$ if and only if $x\in \R\times [-1,1]$.}
\end{example}

Let $\ph: X\to (-\infty, \infty]$ be an extended real-valued function and let $\bar x\in \mbox{dom }\ph$. Recall that $\ph$ is lower semicontinuous at $\bar x$ if
\begin{equation*}
\liminf_{x\to\bar x}\ph(x)\geq \ph(\bar x).
\end{equation*}
The function $\ph$ is called lower semicontinuous if it is lower semicontinuous at every point of its domain.

\begin{theorem} In the same setting of Theorem \ref{t1}, the function $\mathcal{T}^F_Q$ is lower semicontinuous.
\end{theorem}
\noindent{\bf Proof. }Fix any $\bar x\in\mbox{dom }\mathcal{T}^F_Q$ and fix a sequence $\{x_n\}$ that converges to $\bar x$. We will prove that
\begin{equation*}
\liminf \mathcal{T}^F_Q(x_n)\geq \mathcal{T}^F_Q(\bar x)
\end{equation*}
under the assumption (i) of Theorem \ref{t1}. The inequality holds obviously if $\liminf \mathcal{T}^F_Q(x_n)=\infty$. So we only need to consider the case where $\liminf \mathcal{T}^F_Q(x_n)=\gamma\in [0,\infty)$. It suffices to show that $\gamma\geq \mathcal{T}^F_Q(\ox)$. Without loss of generality, we can assume that $\mathcal{T}^F_Q(x_n)\to \gamma$. By the definition of the minimal time function (\ref{MT}), there exists a sequence $\{ t_n\} \subset [0,\infty)$ such that
\begin{equation*}
\mathcal{T}^F_Q(x_n)\leq t_n< \mathcal{T}^F_Q(x_n)+\dfrac{1}{n} \mbox{ and } (x_n+t_n F)\cap  Q\neq \emptyset \mbox{ for all }n\in\N.
\end{equation*}
For every $n\in\N$, fix $f_n\in F$ and $\omega_n\in  Q$ such that
\begin{equation*}
x_n +t_nf_n=\omega_n,
\end{equation*}
and we can assume without loss of generality that $\omega_n\xrightarrow{w} \omega\in  Q$. We will consider two cases: $\gamma>0$ and $\gamma=0$. If $\gamma>0$, then
\begin{equation*}
f_n\xrightarrow{w} \dfrac{\omega-\ox}{\gamma}\in F.
\end{equation*}
Thus $(\ox+\gamma F)\cap  Q\neq \emptyset$, and hence $\gamma\geq \mathcal{T}^F_Q(\ox)$.

Now we consider the case where $\gamma=0$. In this case the sequence $\{t_n\}$ converges to $0$, and
\begin{equation*}
t_nf_n\xrightarrow{w} \omega-\bar x.
\end{equation*}
By Proposition \ref{c2}, $\omega-\bar x\in F_\infty$. This implies $\bar x\in \Omega-F_\infty$. Employing Theorem \ref{t1}, one has that $\mathcal{T}^F_Q(\ox)=0$ and we also have $\gamma\geq \mathcal{T}^F_Q(\ox)$. We have showed that $\mathcal{T}^F_Q$ is lower semicontinuous under the assumption (i) of Theorem \ref{t1}. The proof for the lower semicontinuity under the assumption (ii) is similar.  $\h$ \vspace*{0.05in}

\begin{theorem}\label{mtconvexity} Suppose $ Q \subseteq X$ is a nonempty convex set. Then the minimal time function (\ref{MT}) is convex.
\end{theorem}
\noindent{\bf Proof. }Fix any $x_1, x_2\in \mbox{ dom }\mathcal{T}^F_ Q$ and let $\lambda\in (0,1)$. We will show that
\begin{equation}\label{ci}
\mathcal{T}^F_ Q(\lambda x_1+(1-\lambda)x_2)\leq \lambda \mathcal{T}^F_{ Q}(x_1)+(1-\lambda)\mathcal{T}^F_ Q(x_2).
\end{equation}
Let $\gamma_i=\mathcal{T}^F_ Q(x_i)$ for $i=1,2$. Given any $\epsilon>0$, there exist $t_i$ for $i=1,2$ with
\begin{equation*}
\gamma_i\leq t_i <\gamma_i+\ve \mbox{ and } (x_i+t_iF)\cap  Q\neq \emptyset.
\end{equation*}
Since both $F$ and $ Q$ are convex, one has
\begin{equation*}
[\lambda x_1+(1-\lambda)x_2+(\lambda t_1+(1-\lambda)t_2)F]\cap Q\neq \emptyset.
\end{equation*}
This implies
\begin{equation*}
\mathcal{T}^F_ Q(\lambda x_1+(1-\lambda)x_2)\leq \lambda t_1+(1-\lambda)t_2 \leq \lambda \mathcal{T}^F_ Q(x_1)+(1-\lambda)\mathcal{T}^F_ Q(x_2)+\epsilon.
\end{equation*}
Since $\epsilon$ is arbitrary, (\ref{ci}) holds and the proof is now complete. $\h$ \vspace*{0.05in}

The Minkowski function associated with $F$ is defined by
\begin{equation}\label{mk}
\rho_F(x)=\inf\{ t\geq 0: x\in tF\}, \; x\in X.
\end{equation}
In contrast to the familiar definition of the Minkowski function, in this definition we only assume that $F$ is a nonempty closed convex subset of $X$.

\begin{proposition}\label{hma} The Minkowski function $\rho_F$ is a positively homogenous and subadditive extended real-valued function. Suppose further that $0\in F$. Then  $\rho_F(x)=0$ if and only if $x\in F_\infty$.
\end{proposition}
\noindent{\bf Proof. }We first prove that $\rho_F$ is subadditive, that is,
\begin{equation}\label{add}
\rho_F(x_1+x_2)\leq \rho_F(x_1)+\rho_F(x_2) \mbox{ for all }x_1, x_2\in X.
\end{equation}
The inequality holds obviously if the right side is infinity, that means $x_1 \notin \mbox{dom }\rho_F$ or $x_2\notin\mbox{dom }\rho_F$. Suppose $x_1, x_2\in \mbox{dom }\rho_F$. Fix any $\epsilon>0$. Then there exist $t_1, t_2\geq 0$ such that
\begin{equation*}
\rho_F(x_1)\leq t_1< \rho_F(x_1)+\epsilon, \,\,\, \rho_F(x_2)\leq t_2<\rho_F(x_2)+\epsilon,
\end{equation*}
and
\begin{equation*}
x_1\in t_1F, \; x_2\in t_2F.
\end{equation*}
Since $F$ is convex,
\begin{equation*}
x_1+x_2\in t_1F+t_2F \subseteq (t_1+t_2)F.
\end{equation*}
Thus $\rho_F(x_1+x_2)\leq t_1+t_2<\rho_F(x_1)+\rho_F(x_2)+2\epsilon$. This implies (\ref{add}).

The fact that $\rho_F$ is positive homogeneous follows from the following analysis. For $\alpha > 0$, one has
\begin{align*}
\rho_F(\alpha x)&=\inf\{t\geq 0: \alpha x\in tF\}\\
&=\inf\{t\geq 0: x\in \dfrac{t}{\alpha} F\}\\
&=\alpha \inf\{\dfrac{t}{\alpha} \geq 0: x\in \dfrac{t}{\alpha}F\}=\alpha\rho_F(x).
\end{align*}
Now, we will prove the second statement. Let $ Q=\{0\}$. We have $\rho_F(x)=\mathcal{T}^{-F}_Q(x)$. Since $0\in F$, by Theorem \ref{t1}, $\rho_F(x)=0$ if and only if $x\in  Q-(-F)_\infty=F_\infty$. The proof is now complete. $\h$\vspace*{0.05in}

Recall that a function $\ph: X\to \R$ is said to be $\ell-$Lipschitz if
\begin{equation*}
|\ph(x)-\ph(y)|\leq \ell ||x-y|| \mbox{ for all }x,y\in X.
\end{equation*}
\begin{proposition}\label{lp} Suppose $0\in \mbox{int }F$. Define
\begin{equation*}
\ell=\inf \{ \dfrac{1}{r}: \B(0; r)\subset F, r>0 \}.
\end{equation*}
Then $\rho_F$ is an $\ell-$Lipschitz function. In particular, $\rho_F(x)\leq \ell ||x||$ for all $x\in X$.
\end{proposition}
\noindent{\bf Proof. }It is clear that $\ell$ is a nonnegative real number since $0\in \mbox{int }F$. Let $r>0$ satisfy $\B(0; r)\subset F$. Then
\begin{equation*}
\rho_F(x)=\inf\{ t\geq 0: x\in tF\}\leq \inf\{t\geq 0: x\in t \B(0; r)\}=\inf\{t\geq 0: ||x||\leq rt\}=\dfrac{||x||}{r}.
\end{equation*}
This implies $\rho_F(x)\leq \ell ||x||$. Since $\rho_F$ is subadditive, for any $x, y\in X$, one has
\begin{equation*}
\rho_F(x)-\rho_F(y)\leq \rho_F(x-y)\leq \ell ||x-y||.
\end{equation*}
This implies the $\ell-$Lipschitz property of $\rho_F$. $\h$\vspace*{0.05in}

Using the definitions of the minimal time function (\ref{MT}) and the Minkowski function (\ref{mk}), we can establish a relationship between these functions in the proposition below.

\begin{proposition}\label{mtmk} Let $Q \subseteq X $ be a nonempty closed set.  Then
\begin{equation}\label{rep}
\mathcal{T}^F_Q(x)=\inf\{\rho_F(\omega-x): \omega\in  Q\}.
\end{equation}
\end{proposition}
\noindent{\bf Proof. }Let us first consider the case where $\mathcal{T}^F_Q(x)=\infty$. Then the set
\begin{equation*}
\{ t\geq 0: (x+tF)\cap Q\neq \emptyset\}=\emptyset.
\end{equation*}
It follows that for every $\omega\in Q$, the set
\begin{equation*}
\{t\geq 0: \omega-x\in tF\}=\emptyset.
\end{equation*}
Thus $\rho_F(\omega-x)=\infty$ for all $\omega\in Q$, and hence the right side of (\ref{rep}) is also infinity.

Now, suppose $\mathcal{T}^F_Q(x)<\infty$. Fix any $t\geq 0$ such that
\begin{equation*}
(x+tF)\cap Q\neq \emptyset.
\end{equation*}
Then there exist $f\in F$ and $\omega\in Q$ with $x+tf=\omega$. Thus $\omega-x\in tF$, and hence $\rho_F(\omega-x)\leq t$. This implies
\begin{equation*}
\inf_{\omega\in Q}\rho_F(\omega-x)\leq t.
\end{equation*}
It follows by the definition of the minimal time function (\ref{MT}) that
\begin{equation*}
\inf_{\omega\in Q}\rho_F(\omega-x)\leq \mathcal{T}^F_Q(x)<\infty.
\end{equation*}
Let $\gamma=\inf_{\omega\in Q}\rho_F(\omega-x)$. Then for any $\ve>0$, there exists $\omega\in Q$ with
\begin{equation*}
\rho_F(\omega-x)<\gamma +\ve.
\end{equation*}
Employing the definition of the Minkowski function (\ref{mk}), one finds $t\in [0,\infty)$ with
\begin{equation*}
\omega-x\in tF \mbox{ and }t<\gamma+\ve.
\end{equation*}
Thus $\omega\in x+tF$, and hence
\begin{equation*}
\mathcal{T}^F_Q(x)\leq t <\gamma +\ve.
\end{equation*}
Since $\ve$ is arbitrary, one has
\begin{equation*}
\mathcal{T}^F_Q(x)\leq \gamma=\inf_{\omega\in Q}\rho_F(\omega-x).
\end{equation*}
The proof is now complete. $\h$ \vspace*{0.05in}

Employing the $\ell-$Lipschitz property of the Minkowski function (\ref{mk}) from Proposition \ref{lp} and the representation of the minimal time function (\ref{MT}) in Proposition \ref{mtmk}, we obtain the result below.

\begin{theorem} Let $\ell$ be defined in Proposition \ref{lp}. Suppose $0\in \mbox{int }F$.  Then the minimal time function $\mathcal{T}^F_Q$ is $\ell-$Lipschitz.
\end{theorem}
\noindent{\bf Proof. }Let $x, y\in X$. Then using the subadditivity of $\rho_F$ from Proposition \ref{hma} and using Proposition \ref{lp}, one has
\begin{equation*}
\rho_F(\omega-x)\leq \rho_F(\omega-y)+\rho_F(y-x)\leq \rho_F(\omega-y)+\ell ||y-x|| \mbox{ for all }\omega\in  Q.
\end{equation*}
It follows from Proposition \ref{mtmk} that
\begin{equation*}
\mathcal{T}^F_Q(x)\leq \mathcal{T}^F_Q(y)+\ell ||x-y||.
\end{equation*}
This implies the $\ell-$Lipschitz property of $\mathcal{T}^F_Q$. $\h$

\section{Convex Subdifferentials of Minimal Time Functions with Unbounded Dynamics}
\label{sec:3}

In this section, we discuss properties of the convex subdifferentials of the minimal time function (\ref{MT}) where the target set $Q$ is convex.  In this case the minimal time function (\ref{MT}) is also convex by Theorem \ref{mtconvexity}.
Consider the following extended real-valued function defined on $X^*$:
\begin{equation*}
\sigma_F(x^*)=\sup\{\la x^*, f\ra: f\in F\}.
\end{equation*}

\begin{theorem}\label{inset1} Consider the minimal time function (\ref{MT}) with a nonempty convex target set $Q$. For $\bar x\in  Q$, one has
\begin{equation*}
\partial \mathcal{T}^F_Q(\ox)=N(\ox;  Q)\cap C^*,
\end{equation*}
where
\begin{equation*}
C^*=\{x^*\in X^*: \sigma_F(-x^*)\leq 1\}.
\end{equation*}
Suppose further that $0\in F$. Then
\begin{equation*}
\partial \mathcal{T}^F_Q(\ox)=N(\ox; Q-F_\infty)\cap C^*=N(\ox;  Q)\cap C^*.
\end{equation*}
\end{theorem}
\noindent{\bf Proof. }Fix any $x^*\in \partial \mathcal{T}^F_ Q(\ox)$. One has
\begin{equation}\label{convexity}
\la x^*, x-\bar x\ra \leq \mathcal{T}^F_ Q(x)-\mathcal{T}^F_ Q(\ox)\mbox{ for all }x\in X.
\end{equation}
Since $\mathcal{T}^F_ Q(x)=0$ for all $x\in Q$, one has
\begin{equation*}
\la x^*, x-\ox\ra \leq 0\mbox{ for all }x\in  Q.
\end{equation*}
This implies $x^*\in N(\ox;  Q)$. For any $f\in F$ and $t>0$, it follows from (\ref{convexity}) that
\begin{equation*}
\la  x^*, (\ox-tf)-\ox\ra \leq \mathcal{T}^F_Q(\ox-tf)\leq t.
\end{equation*}
Note that the last inequality holds since $((\ox-tf)+tF)\cap  Q\neq \emptyset$. Thus
\begin{equation*}
\la x^*, -f\ra \leq 1 \mbox{ for all }f\in F,
\end{equation*}
so $x^*\in C^*$.  Thus $\partial \mathcal{T}^F_Q(\ox) \subset N(\ox;  Q)\cap C^*$.

Now, we will show that
\begin{equation*}
N(\ox;  Q)\cap C^* \subset \partial \mathcal{T}^F_ Q(\ox).
\end{equation*}
Fix any $x^*\in N(\ox;  Q)\cap C^*$. Then $\la x^*, -f\ra \leq 1$ for all $f\in F$ and
\begin{equation*}
\la x^*, x-\ox\ra \leq 0\mbox{ for all }x\in  Q.
\end{equation*}
Fix any $u\in \mbox{ dom }\mathcal{T}^F_Q$. For any $\ve>0$, there exist $t\in [0,\infty)$, $f\in F$, and $\omega\in  Q$ with
\begin{equation*}
\mathcal{T}^F_Q(u)\leq t< \mathcal{T}^F_Q(u)+\ve \mbox{ and } u+tf=\omega.
\end{equation*}
Then
\begin{align*}
\la x^*, u-\ox\ra &=\la x^*, u-\omega\ra +\la x^*, \omega-\ox\ra\\
&\leq \la x^*, -tf\ra \leq t< \mathcal{T}^F_Q(u)+\ve =\mathcal{T}^F_Q(u)-\mathcal{T}^F_ Q(\ox)+\ve.
\end{align*}
This implies (\ref{convexity}) since $\ve$ is arbitrary.

Let us finally prove that
\begin{equation*}
\partial \mathcal{T}^F_ Q(\ox)=N(\ox; Q-F_\infty)\cap C^*=N(\ox;  Q)\cap C^*,
\end{equation*}
under the additional assumption. Obviously, $N(\ox; Q-F_\infty)\cap C^*\subset N(\ox;  Q)\cap C^*$ since $ Q\subset  Q-F_\infty$. Fix $x^*\in N(\ox; Q)\cap C^*$ and any $x\in  Q-F_\infty$. One has
$x=\omega-d$ for $\omega\in Q$ and $d\in F_\infty$. Since $0\in F$, one has $td\in F$ for all $t>0$. This implies $\la x^*, -td\ra\leq 1$ for all $t>0$. Thus $\la x^*, -d\ra \leq 0$. Then
\begin{equation*}
\la x^*, x-\ox\ra =\la x^*, \omega-d-\ox\ra =\la x^*, \omega-\ox\ra +\la x^*, -d\ra \leq 0.
\end{equation*}
Therefore, $x^*\in N(\ox;  Q-F_\infty)\cap C^*$. The proof is now complete. $\h$ \vspace*{0.05in}

Let us now consider the case where $\ox\in  Q-F_\infty$. In fact, similar results like those in Theorem \ref{inset1} hold true in this case.
\begin{theorem}\label{outset1} Consider the minimal time function (\ref{MT}) with a nonempty convex target set $Q$. Suppose $\ox\in  Q-F_\infty$ and $0\in F$. Then
\begin{equation*}
\partial \mathcal{T}^F_ Q(\ox)=N(\ox;  Q-F_\infty)\cap C^*.
\end{equation*}
\end{theorem}
\noindent{\bf Proof. }Fix any $x^*\in \partial \mathcal{T}^F_ Q(\ox)$. Using (\ref{convexity}) and the fact that $\mathcal{T}^F_ Q(x)=0$ for all $x\in  Q-F_\infty$ from Remark \ref{r1}, one has that $x^*\in N(\ox; Q-F_\infty)$. Since $\ox\in  Q-F_\infty$, one has $\ox=\bar \omega-d$, where $\bar\omega\in Q$ and $d\in F_\infty$. Fix any $f\in F$. For any $t>0$, let $x_t=(\ox+d)-tf=\bar\omega-tf$. Then $(x_t+tF)\cap Q\neq \emptyset$. Thus
\begin{equation*}
\la x^*, x_t-\ox\ra \leq \mathcal{T}^F_Q(x_t)\leq t,
\end{equation*}
or equivalently,
\begin{equation*}
\la x^*, d-tf\ra\leq  t.
\end{equation*}
Thus
\begin{equation*}
\la x^*, \dfrac{d}{t}-f\ra\leq  1 \mbox{ for all }t>0.
\end{equation*}
This implies $\la x^*, -f\ra\leq 1$ by letting $t\to\infty$, and hence $x^*\in C^*$.
Thus we have shown that $\partial \mathcal{T}^F_ Q(\ox) \subset N(\ox;  Q-F_\infty)\cap C^*$.

Now fix any $x^*\in N(\ox; Q-F_\infty)\cap C^*$. Fix any $u\in \mbox{dom }\mathcal{T}^F_Q$. Then there exist $t\in [0,\infty)$, $\omega\in  Q$ and $f\in F$ such that
\begin{equation*}
\mathcal{T}^F_Q(u)\leq t<\mathcal{T}^F_Q(u)+\ve \mbox{ and } u+tf=\omega.
\end{equation*}
Since $ Q\subset  Q-F_\infty$, one has
\begin{equation*}
\la x^*, u-\ox\ra =\la x^*, \omega-\ox\ra +t\la x^*, -f\ra \leq t<\mathcal{T}^F_Q(u)+\ve=\mathcal{T}^F_Q(u)-\mathcal{T}^F_ Q(\ox)+\ve.
\end{equation*}
This again implies $x^*\in \partial \mathcal{T}^F_ Q(\ox)$, and thus $N(\ox;  Q-F_\infty)\cap C^* \subset  \partial \mathcal{T}^F_ Q(\ox)$. The proof is now complete. $\h$ \vspace*{0.05in}

In what follows we are going to discuss the subdifferential in the
sense of convex analysis for the minimal time function (\ref{MT})
when the reference point does not belong to $Q-F_\infty$. Our
development is inspired by the early work on subdifferentials of
distance functions; see \cite{boun-th02,BJQ,K}. For any $r>0$,
define
\begin{equation*}
 Q_r:=\{x\in X: \mathcal{T}^F_ Q(x)\leq r\}.
\end{equation*}
\begin{lemma}\label{distanceestimate} Let $Q$ be a nonempty set and let $r>0$. Suppose $x\notin  Q_r$ with $\mathcal{T}^F_ Q(x)<\infty$. Then
\begin{equation*}
\mathcal{T}^F_ Q(x)=\mathcal{T}^F_{Q_r}(x)+r.
\end{equation*}
\end{lemma}
\noindent{\bf Proof. }Let $r >0$. Since $Q\subset Q_r$, one has $\mathcal{T}^F_{Q_r}(x)\leq \mathcal{T}^F_Q(x)<\infty$. Fix any $t\geq 0$ such that $$(x+tF)\cap Q_r\neq\emptyset.$$
Then there exist $f_1\in F$ and $u\in Q_r$ such that $x+tf_1=u$. Since $u \in Q_r$, then $\mathcal{T}^F_Q(u)\leq r$. Then for any $\ve>0$, there exists $s>0$ such that $s<r+\ve$. Consequently, there exist $\omega\in Q$ and $f_2\in F$ such that $u+sf_2=\omega$. Thus
\begin{equation*}
\omega=u+sf_2 = (x+tf_1)+sf_2\in x+(t+s)F
\end{equation*}
since $F$ is convex.  This implies
\begin{equation*}
\mathcal{T}^F_Q(x)\leq t+s\leq t+r+\ve.
\end{equation*}
Since $\ve>0$ is arbitrary, $\mathcal{T}^F_Q(x)\leq t+r$, and hence $$\mathcal{T}^F_Q(x)\leq \mathcal{T}^F_{Q_r}+r.$$

Let us prove the opposite inequality. Let $\gamma=\mathcal{T}^F_Q(x)$. Then $r<\gamma$ since $x\notin Q_r$. For any $\ve>0$, there exist $t\in [0,\infty)$, $f\in F$, and $\omega\in Q$ with
\begin{equation*}
\gamma\leq t<\gamma+\ve \mbox{ and } x+tf=\omega.
\end{equation*}
One has
\begin{equation*}
\omega=x+tf=x+(t-r)f+rf \in x+(t-r)f +rF.
\end{equation*}
Thus
\begin{equation*}
\mathcal{T}^F_Q(x+(t-r)f)\leq r.
\end{equation*}
So $x+(t-r)f\in Q_r$; in addition $x+(t-r)f\in x + (t-r)F$. Hence $\mathcal{T}^F_{Q_r}(x)\leq t-r\leq \gamma-r+\ve$. Since $\ve>0$ is arbitrary, one has
\begin{equation*}
r+\mathcal{T}^F_{Q_r}\leq \gamma =\mathcal{T}^F_Q(x).
\end{equation*}
The proof is now complete. $\h$

\begin{lemma}\label{lm3} Let $Q$ be a nonempty set and let $x\in\mbox{dom }\mathcal{T}^F_Q$. For $ t \geq 0$ and any $f\in F$, one has
\begin{equation*}
\mathcal{T}^F_Q(x-tf)\leq \mathcal{T}^F_ Q(x)+t.
\end{equation*}
\end{lemma}
\noindent{\bf Proof. }For any $\ve>0$, there exists $s \geq 0$ such that
\begin{equation*}
\mathcal{T}^F_ Q(x)\leq s<\mathcal{T}^F_ Q(x)+\ve \mbox{ and } (x+sF)\cap  Q\neq \emptyset.
\end{equation*}
Then $(x-tf + tF+sF)\cap  Q\neq \emptyset$, and hence $(x-tf + (t+s)F)\cap  Q\neq \emptyset$. Thus
\begin{equation*}
\mathcal{T}^F_Q(x-tf)\leq t+s\leq \mathcal{T}^F_Q(x)+t+\ve.
\end{equation*}
The conclusion follows by letting $\ve\to 0$. $\h$

\begin{theorem}\label{outset} Let $X$ be reflexive Banach space and let $0\in F$. Consider the minimal time function (\ref{MT}) with a nonempty closed bounded convex target set $Q$. For $\bar x\notin  Q-F_\infty$, one has
\begin{equation*}
\partial \mathcal{T}^F_Q(\ox)=N(\ox;  Q_r)\cap S^*,
\end{equation*}
where
\begin{equation*}
S^*=\{x^*\in X^*: \sigma_F(-x^*)= 1\} \mbox{ and }r=\mathcal{T}^F_ Q(\ox)>0.
\end{equation*}
\end{theorem}
 \noindent{\bf Proof. }Fix any $x^*\in \partial \mathcal{T}^F_Q(\ox)$. Following the proof of Theorem \ref{inset1}, see that $\sigma_F(-x^*)\leq 1$ and $x^*\in N(\ox;  Q_r)$. Let us show that $\sigma_F(-x^*)=1$. One has
 \begin{equation}\label{convexity1}
 \la x^*, x-\ox\ra \leq \mathcal{T}^F_ Q(x)-\mathcal{T}^F_ Q(\ox) \mbox{ for all }x\in X.
 \end{equation}
 Fix any $\ve\in (0, r)$. There exist $t\in \R$, $f\in F$ and $\omega\in  Q$ such that
 \begin{equation*}
 r\leq t<r+\ve^2 \mbox{ and }\omega=\ox+tf.
 \end{equation*}
 We can write $\omega=\ox+\ve f +(t-\ve)f$. So $\mathcal{T}^F_Q(\ox+\ve f)\leq t-\ve$. Applying (\ref{convexity1}) for $x=\ox+\ve f$, one has
 \begin{equation*}
 \la x^*, \ve f\ra \leq \mathcal{T}^F_Q(\ox+\ve f)-\mathcal{T}^F_ Q(\ox)\leq t-\ve -r\leq \ve^2-\ve.
 \end{equation*}
 This implies
 \begin{equation*}
 1-\ve \leq \la -x^*, f\ra \leq \sigma_F(-x^*).
 \end{equation*}
 So $\sigma_F(-x^*)\geq 1$ by letting $\ve\to 0$. Therefore, $x^*\in S^*$, and thus $\partial\mathcal{T}^F_Q(\ox)\subset N(\ox;Q_r)\cap S^*$.

 Fix any $x^*\in N(\ox; Q_r)$ such that $\sigma_F(-x^*)=1$. We are going to show that (\ref{convexity1}) is satisfied. By Theorem \ref{inset1}, $x^*\in \partial \mathcal{T}^F_{Q_r}(\ox)$. Thus
 \begin{equation*}
 \la x^*, x-\ox\ra \leq \mathcal{T}^F_{Q_r}(x) \mbox{ for all }x\in X.
 \end{equation*}
 Fix any $x\in X$. In the case $t:=\mathcal{T}^F_ Q(x)>r$, one has $\mathcal{T}^F_ Q(x)-r=\mathcal{T}^F_{Q_r}(x)$ by Lemma \ref{distanceestimate}. Thus (\ref{convexity1}) holds. Suppose $t\leq r$. For any $\ve>0$, choose $f\in F$ such that $\la x^*, -f\ra >1-\ve$. By Lemma \ref{lm3}, $\mathcal{T}^F_Q(x-(r-t)f)\leq r$. So $x-(r-t)f\in  Q_r$. Since $x^*\in N(\ox; Q_r)$, one has
 \begin{equation*}
 \la x^*, x-(r-t)f-\ox\ra \leq 0.
 \end{equation*}
 This implies
 \begin{equation*}
 \la x^*, x-\ox\ra\leq \la x^*, f\ra (r-t)\leq (1-\ve)(t-r)=(1-\ve)(\mathcal{T}^F_ Q(x)-\mathcal{T}^F_ Q(\ox)).
 \end{equation*}
 Since $\ve$ is arbitrary, (\ref{convexity1}) holds. Therefore, $x^*\in \partial \mathcal{T}^F_ Q(\ox)$. $\h$

 \begin{theorem}\label{subofs} Let $X$ be reflexive Banach space and let $0\in F$. Consider the minimal time function (\ref{MT}) with a nonempty closed bounded convex target set $Q$. Suppose $\ox\notin  Q-F_\infty$ and $0\in F$. Then
 \begin{equation*}
 \partial \mathcal{T}^F_ Q(\ox) =[-\partial \rho_F(\bar\omega-\ox)]\cap N(\bar\omega; Q)
 \end{equation*}
 for any $\bar\omega\in \Pi_F(\ox; Q)$, where
 \begin{equation*}
 \Pi_F(\ox; Q)=\{\omega\in  Q: \mathcal{T}^F_ Q(\ox)=\rho_F(\omega-\ox)\}.
 \end{equation*}
 \end{theorem}
 \noindent{\bf Proof. }Fix any $x^*\in \partial \mathcal{T}^F_ Q(\ox)$ and $\bar\omega\in \Pi_F(\ox;Q)$. Then $\mathcal{T}^F_Q(\ox)=\rho_F(\bar\omega-\ox)$ and
 \begin{equation}\label{subdifferentialinequality}
 \la x^*, x-\ox\ra \leq \mathcal{T}^F_Q(x)-\mathcal{T}^F_Q(\ox) \mbox{ for all }x\in X.
  \end{equation}
  For any $\omega\in Q$, one has
  \begin{align*}
  \la x^*, \omega-\bar\omega\ra& =\la x^*, (\omega-\bar\omega+\ox)-\ox\ra\\
 & \leq \mathcal{T}^F_Q(\omega-\bar\omega+\ox)-\mathcal{T}^F_Q(\ox)\\
  &\leq \rho_F(\omega-(\omega-\bar\omega+\ox))-\rho_F(\bar\omega-\ox)=0.
  \end{align*}
  This implies $x^*\in N(\bar\omega;Q)$. Let $\Tilde{u}=\bar\omega-\ox$. For any $t\in (0,1)$ and $u\in X$, applying (\ref{subdifferentialinequality}) with $x=\ox-t(u-\Tilde{u})$, one has
\begin{align*}
\la x^*, -t(u-\Tilde{u})\ra &\leq \mathcal{T}^F_Q(\ox-t(u-\Tilde{u}))-\mathcal{T}^F_Q(\ox)\\
&\leq\rho_F(\bar\omega-(\ox-t(u-\Tilde{u})))-\rho_F(\bar\omega-\ox)\\
&=\rho_F(\Tilde{u}+t(u-\Tilde{u}))-\rho_F(\Tilde{u})\\
&=\rho_F(tu+(1-t)\Tilde{u})-\rho_F(\Tilde{u})\\
&\leq t\rho_F(u)+(1-t)\rho_F(\Tilde{u})-\rho_F(\Tilde{u})\\
&=t(\rho_F(u)-\rho_F(\Tilde{u})).
\end{align*}
It follows that
\begin{equation*}
\la -x^*, u-\Tilde{u}\ra\leq \rho_F(u)-\rho_F(\Tilde{u})\mbox{ for all }u\in X.
\end{equation*}
Thus $-x^*\in \partial \rho_F(\Tilde{u})$, and hence $x^*\in -\partial\rho_F(\bar\omega-\ox)$. The inclusion ``$\subset$" has been proved. Let us prove the opposite inclusion. We have that
 \begin{equation*}
 \rho_F(x)=\inf\{t\geq 0: x\in tF\}=\inf\{t\geq 0: (-x+tF)\cap O\neq \emptyset\}=\mathcal{T}^F_O(-x),
 \end{equation*}
 where $O=\{0\}$. For any $\bar\omega\in Q$ one has, $\ox-\bar\omega\notin O-F_\infty=-F_\infty$ since $\ox\notin Q-F_\infty$. From $-x^*\in \partial \rho_F(\bar\omega-\ox)$, one has $x^*\in S^*$ by Theorem \ref{outset}. Thus we only need to show that
 \begin{equation*}
 [-\partial \rho_F(\bar\omega-\ox)]\cap N(\bar\omega; Q)\subset N(\ox; Q_r),
 \end{equation*}
 and the conclusion will follow from Theorem \ref{outset}.

To proceed, pick any $x\in  Q_r$ and for an arbitrary small $\ve>0$. Then there exist $t<r+\ve$, $f\in F$,
and $\omega\in Q$ with $\omega=x+tf$. Thus $\la-x^*,f\ra\le\sigma_F(-x^*)= 1$ and
\begin{align*}
\la x^*,x-\ox\ra&=\la x^*,\omega-tf-\ox\ra\\
&=t\la-x^*,f\ra+\la x^*,\omega-\bar\omega\ra +\la x^*,
\bar\omega-\ox\ra\\
&\le t+\la x^*, \omega-\bar\omega\ra+\la x^*,\bar\omega-\ox\ra\\
&\le \mathcal{T}^F_ Q(\ox)+\ve+\la x^*,\omega-\bar\omega\ra+\la x^*,
\bar\omega-\ox\ra.
\end{align*}
We have $\la x^*,\omega-\bar\omega\ra\le 0$ due to $x^*\in N(\bar\omega; Q)$ and
\begin{align*}
\la x^*,\bar\omega-\ox\ra=\la-x^*,0-(\bar\omega-\ox)\ra\le
\rho_F(0)-\rho_F(\bar\omega-\ox)=-\mathcal{T}^F_ Q(\ox)
\end{align*}
by $-x^*\in \partial\rho_F(\bar\omega-\ox)$. It follows therefore that $\la x^*,x-\ox\ra\le\ve$ for all $x\in Q_r$,
and hence $x^*\in N(\ox; Q_r)$ because $\ve>0$ was chosen arbitrary small.  $\h$\vspace*{0.05in}

Let $\ph: X\to (-\infty, \infty]$ be a convex function and let $\ox\in \mbox{dom }\ph$. The singular subdifferential of $\ph$ at $\ox$ is defined by
\begin{equation*}
\partial^\infty \ph(\ox)=\{x^*\in X^*: (x^*, 0)\in N((\ox, \ph(\ox)); \epi \ph)\}.
\end{equation*}
\begin{lemma}\label{singular} Let $\ph: X\to (-\infty, \infty]$ be a convex function and let $\ox\in \mbox{dom }\ph$. We always have
\begin{equation*}
\partial^\infty \ph(\ox)=N(\ox; \mbox{dom }\ph).
\end{equation*}
\end{lemma}
\noindent{\bf Proof. }Fix any $x^*\in \partial^\infty \ph(\ox)$ and $x\in \mbox{dom }\ph$. Then $(x, \ph(x))\in \epi \ph$. Thus
\begin{equation*}
\la x^*, x-\ox\ra =\la x^*, x-\ox\ra +0(\ph(x)-\ph(\ox))\leq 0
\end{equation*}
because $(x^*, 0)\in N((\ox, \ph(\ox)); \epi \ph)$. Conversely, suppose $x^*\in N(\ox; \mbox{dom }\ph)$. Fix any $(x, \lambda)\in \epi f$. Then $\ph(x)\leq \lambda$, so $x\in \mbox{dom }\ph$. Thus
\begin{equation*}
\la x^*, x-\ox\ra +0(\lambda-\ph(\ox))=\la x^*, x-\ox\ra \leq 0.
\end{equation*}
This implies $(x^*, 0)\in N((\ox, \ph(\ox)); \epi \ph)$ or equivalently $x^*\in\partial^\infty \ph(\ox)$. $\h$ \vspace*{0.05in}

Define the following subset of $X^*$
\begin{equation*}
F^*_{+}=\{x^*\in X^*: \la x^*, f\ra \geq 0 \mbox{ for all }f\in F\}.
\end{equation*}
We are going to obtain below explicit representations of singular subdifferentials of the minimal time function (\ref{MT}).

\begin{proposition} Consider the minimal time function (\ref{MT}) with a nonempty convex target set $ Q$. The following hold:\\
{\rm (i)} If $\ox\in  Q$, then
\begin{equation*}
\partial^\infty \mathcal{T}^F_Q(\ox)= N(\ox;  Q)\cap F^*_{+}.
\end{equation*}
{\rm (ii)} If $\ox\in  Q-F_\infty$ and $0\in F$, then
\begin{equation*}
\partial^\infty \mathcal{T}^F_ Q(\ox)=N(\ox;  Q-F_\infty)\cap F^*_{+}.
\end{equation*}
{\rm (iii)} If $\ox\notin  Q-F_\infty$, $X$ is a reflexive Banach space, and $Q$ is closed and bounded, $0\in F$, and $\mathcal{T}^F_ Q(\ox)<\infty$, then
\begin{equation*}
\partial^\infty \mathcal{T}^F_ Q(\ox)=N(\ox;  Q_r)\cap F^*_{+},\; r=\mathcal{T}^F_ Q(\ox)>0.
\end{equation*}
\end{proposition}
\noindent{\bf Proof.} (i) Fix any $x^*\in \partial^\infty \mathcal{T}^F_Q(\ox)=N(\ox;\mbox{dom }\mathcal{T}^F_Q)$. Employing Lemma \ref{singular} and the fact that $ Q\subset \mbox{dom }\mathcal{T}^F_Q$, one has $\partial^\infty \mathcal{T}^F_Q(\ox)=N(\ox;\mbox{dom }\mathcal{T}^F_Q)\subset N(\ox; Q)$. We are going to show that $\la x^*, f\ra \geq 0$ for all $f\in F$. For any $f\in F$, one has $\ox -f \in \mbox{dom }\mathcal{T}^F_Q$ since $\mathcal{T}^F_Q(\ox-f)\leq 1$. Thus
\begin{equation*}
\la x^*, (\ox-f)-\ox\ra \leq 0.
\end{equation*}
This implies $\la x^*, f\ra \geq 0$. Let us prove the converse. Fix any $x^*\in N(\ox;  Q)\cap F^*_{+}$ and $x\in \mbox{dom }\mathcal{T}^F_Q$. Then there exists $t\in [0,\infty)$ such that
$(x+tF)\cap Q\ne\emp$. Fix further $f\in F$ and $\omega\in Q$ with $x+tf=w$. Then
\begin{align*}
\la x^*,x-\ox\ra&=\la x^*,\omega-tf-\ox\ra\\
&=\la x^*,\omega-\ox\ra-t\la x^*,f\ra\le 0,
\end{align*}
Thus we get $x^*\in N(\ox; \mbox{dom }\mathcal{T}^F_Q)=\partial^\infty \mathcal{T}^F_Q(\ox)$.\\[1ex]
(ii) Since $\mathcal{T}^F_ Q(x)=0$ for all $x\in  Q-F_\infty$, one has $ Q-F_\infty\subset \mbox{dom }\mathcal{T}^F_Q$. Applying Lemma \ref{singular} again, one has $\partial^\infty\mathcal{T}^F_Q(\ox)=N(\ox; \mbox{dom }\mathcal{T}^F_Q)\subseteq N(\ox;  Q-F_\infty)$. Let $\bar\omega\in  Q$ and $d\in F_\infty$ be satisfy $\ox=\bar\omega-d$. Since $0\in F$, one has $d\in F$. For any $f\in F$,
\begin{equation*}
\ox-f =\bar\omega-(d+f)\subset \bar\omega -2(1/2d+1/2f)\subset \bar\omega-2F.
\end{equation*}
This implies $\mathcal{T}^F_Q(\ox-f)\leq 2$, so $\ox-f\in \mbox{dom }\mathcal{T}^F_Q$. This again implies $x^*\in F_{+}$ as in the proof of (i). The proof of the oppositive inclusion follows by the same proof of (i) using the observation that $ Q\subset  Q-F_\infty$. The proof for (iii) is similar. $\h$

\section{Applications to the Smallest Intersecting Ball Problem}
\label{sec:4}

In this section, we consider the optimization problem (\ref{SIB}),
where the target sets $\Omega_i$ for $i=1,\ldots,m$ and the
constraint set $\Omega$ are nonempty closed bounded convex subsets
of $X$. The smallest intersecting ball problem (SIB) under
consideration asks for a point $\ox\in \Omega$ and the smallest
$t\geq 0$ (if they exist) such that
$$(\ox+tF)\cap \Omega_i\neq \emptyset \mbox{ for all } i=1,\ldots,m.$$
\begin{definition}\label{solution} We say that $\ox\in\Omega$ is a solution of the smallest intersecting ball problem (SIB) with the smallest radius $r$ if $r\in [0,\infty)$,
\begin{equation*}
(\ox+rF)\cap\Omega_i\neq \emptyset \mbox{ for all }i=1,\ldots,m,
\end{equation*}
and whenever there exist $x\in \Omega$ and $r^\prime\in [0,\infty)$
with
\begin{equation*}
(x+r^\prime F)\cap\Omega_i\neq \emptyset \mbox{ for all
}i=1,\ldots,m,
\end{equation*}
then $r\leq r^\prime$.
\end{definition}

\begin{lemma}\label{lm} Let $X$ be a reflexive Banach space. Suppose $Q \subseteq X$ is a nonempty closed bounded convex set and $0\in F$. If $\mathcal{T}^F_Q(\ox)\leq r$ and $r>0$, then $(\ox+rF)\cap Q\neq \emptyset$.
\end{lemma}
\noindent{\bf Proof. }Let $\gamma=\mathcal{T}^F_Q(\ox)$. Suppose first that $\gamma>0$.
Then there exists a sequence $\{ t_k\}$ that converges to $\gamma$ such that $\ox+t_kf_k=\omega_k$, where $f_k\in
F$ and $\omega_k\in Q$. Since $Q$ is closed convex and bounded, we can assume that
$\omega_k\xrightarrow{w}\omega\in Q$. Then $f_k\xrightarrow{w}
\dfrac{\omega-\ox}{\gamma}$. Since $F$ is closed and convex, it is weakly closed, and hence
$\dfrac{\omega-\ox}{\gamma}\in F$. Thus $(\ox+\gamma
F)\cap Q\neq \emptyset.$ It follows that $(\ox+r F)\cap Q\neq
\emptyset$ since $\gamma\leq r$ and $0\in F$. If $\gamma=0$, then $\ox\in Q-F_\infty$ by Theorem \ref{t1}. Let
$\ox=\omega-d$, where $\omega\in Q$ and $d\in F_\infty$. By Corollary \ref{c1}, for $r \geq 0$
\begin{equation*}
\dfrac{\omega-\ox}{r}=\dfrac{d}{r}\in F.
\end{equation*}
This again implies $(\ox+rF)\cap Q\neq\emptyset$. The proof is now complete. $\h$

\begin{theorem} Let $X$ be a reflexive Banach space and let $0\in \mbox{int }F$. Suppose $\Omega$ and $\Omega_i$ for $i=1,\ldots,m$ are nonempty closed bounded convex sets and
\begin{equation}\label{emptyinter}
\Omega\cap[\cap_{i=1}^m(\Omega_i-F_\infty)]=\emptyset.
\end{equation}
Then $\ox\in \Omega$ is an optimal solution of the optimization
problem (\ref{SIB}) with $r=\mathcal{T}(\ox)$ if and only if $\ox$
is an optimal solution of the smallest intersecting ball problem
(SIB) with the smallest radius $r$.
\end{theorem}
\noindent{\bf Proof. }Suppose $\ox\in\Omega$ is an optimal solution of
(\ref{SIB}) and $r=\mathcal{T}(\ox)$. If $r=0$, then
$\mathcal{T}^F_{\Omega_i}(\ox)=0$ for all $i=1,\ldots,m$. This
implies $\ox\in \Omega_i-F_\infty$ for all $i=1,\ldots,m$, which is
a contradiction to (\ref{emptyinter}). Thus $r>0$. From Lemma
\ref{lm}, $(\ox+rF)\cap \Omega_i\neq\emptyset$ for all $i=1,\ldots,
m$. Let us fix an $x\in \Omega$ and $r^\prime\geq 0$ such that
$(x+r^\prime F)\cap \Omega_i\neq \emptyset$ for all $i=1,\ldots,m$.
We will show that $r^\prime\geq r$. Suppose $r^\prime<r$. Since
$\mathcal{T}^F_{\Omega_i}(x)\leq r^\prime$ for all $i=1,\ldots,m$,
one has $\mathcal{T}(x)\leq r^\prime <r=\mathcal{T}(\ox)$, which is
a contradiction.

Let us prove the converse. Suppose that $\ox$ is an optimal solution
of the smallest intersecting ball problem (SIB) with the smallest
radius $r$. From Definition \ref{solution}, one has $\ox\in \Omega$
and $(\ox+rF)\cap \Omega_i\neq \emptyset$ for all $i=1,\ldots,m$.
Thus
\begin{equation}\label{1}
\mathcal{T}(\ox)=\max\{\mathcal{T}^F_{\Omega_i}(\ox):
i=1,\ldots,m\}\leq r.
\end{equation}
By (\ref{emptyinter}), $r>0$. Fix any $x\in \Omega$ and define
$r^\prime=\mathcal{T}(x)$. Since $0\in \mbox{int }F$ and by
(\ref{emptyinter}), one has $r^\prime\in (0,\infty)$. By the
definition of the minimal time function (\ref{MT}),
$\mathcal{T}^F_{\Omega_i}(x)\leq r^\prime$ for all $i=1,\ldots,m$.
It follows again from Lemma \ref{lm} that
\begin{equation*}
(x+r^\prime F)\cap \Omega_i\neq \emptyset \mbox{ for all
}i=1,\ldots,m.
\end{equation*}
Thus $r\leq r^\prime=\mathcal{T}(x)$ by Definition \ref{solution}. In particular, taking $x=\ox$, one has $r\leq \mathcal{T}(\ox)$, and hence $r=\mathcal{T}(\ox)$ by (\ref{1}). So $\mathcal{T}(\ox)\leq \mathcal{T}(x)$. Therefore, $\ox$ in an optimal solution of the optimization problem (\ref{SIB}). $\h$ \vspace*{0.05in}

For each $u\in X$, define
\begin{equation*}
I(u)=\{i\in \{1,\ldots,m\}:
\mathcal{T}(u)=\mathcal{T}^F_{\Omega_i}(u)\}.
\end{equation*}

\begin{theorem}\label{algorithm} Let $X=\R^n$ with the Euclidean norm and let $0\in \mbox{int }F$. Consider the optimization problem (\ref{SIB}) in which $\Omega_i$ for $i=1,\ldots,m$ and $\Omega$ are nonempty closed bounded convex sets. Fix $x_1\in \Omega$ and define the sequences of iterates by
\begin{equation*}
x_{k+1}=\mathcal{P}(x_k-\alpha_k x^*_k;\Omega), \, \, \, \forall k
\in \N,
\end{equation*}
where
\begin{eqnarray*}
x^*_k\in\left\{\begin{array}{lr} \{0\}&\mbox{ if
}\;x_k\in\Omega_i-F_{\infty},\\\\

[-\partial \rho_F(\omega_k-x_k)]\cap N(\omega_k;\Omega_i) &\mbox{ if
}\;x_k\notin\Omega_i-F_{\infty},
\end{array}
\right.
\end{eqnarray*}
for some $i\in I(x_k)$, and $\mathcal{P}$ denotes the Euclidean projection operator.

Let $V_k=\min\{\mathcal{T}(x_j): j=1,\ldots,k\}$. If $\sum_{i=1}^\infty
\alpha_i=\infty$ and $\sum_{i=1}^\infty \alpha_i^2<\infty$, then the sequence
$\{x_k\}$ converges to an optimal solution of problem (\ref{SIB}) and
$V_k$ converges to the optimal value $\Bar V$ of the problem.
Moreover,

\begin{align*}
V_k-\Bar
V\le\dfrac{d(x_1;S)^2+L^2\sum_{i=1}^k\alpha_i^2}{2\sum_{i=1}^k\alpha_i},
\end{align*}
where $0\le L<\infty$ is a Lipschitz constant of the function
$\mathcal{T},$ and $S$ is the solution set of problem (\ref{SIB}).

\end{theorem}
\noindent{\bf Proof. }First observe that $S\neq \emptyset$ since
$\mathcal{T}$ is a continuous function and the constraint $\Omega$
is compact. We have the following well-known formula, see e.g.
\cite[Corollary~4.3.2]{hl}, for computing the subdifferential of the
function $\mathcal{T}$ at a point $u$:
\begin{equation*}
\partial \mathcal{T}(u)=\co\{ \partial \mathcal{T}^F_{\Omega_i}(u): i\in I(u)\}.
\end{equation*}
Thus for any $i\in I(u)$, one has $\partial
\mathcal{T}^F_{\Omega_i}(u)\subset \partial\mathcal{T}(u)$. The
choice of $x^*_k$ in the statement of the theorem allows us to find
an element of $\partial \mathcal{T}^F_{\Omega_i}(x_k)$ by Theorem
\ref{inset1}, Theorem \ref{outset1}, and Theorem \ref{subofs}. Then
$x^*_k\in \partial \mathcal{T}(x_k)$. The  rest of the proof comes
from the standard projective subgradient method; see \cite{bert}.
$\h$\vspace*{0.05in}

In the proposition below, we provide explicit formulas for computing the asymptotic cone and the Minkowski function generated by an unbounded dynamic.

\begin{proposition}\label{poly} Let $x^*_i\in X^*$ and $b_i> 0$ for $i=1,\ldots,q$. Define
\begin{equation}\label{pdy}
F=\{x\in X: \la x^*_i, x\ra \leq b_i \mbox{ for all }i=1,\ldots,
q\}.
\end{equation}
Then
\begin{equation*}
F_\infty=\{x\in X: \la x^*_i, x\ra \leq 0 \mbox{ for all
}i=1,\ldots, q\}
\end{equation*}
and
\begin{eqnarray*}
\rho_F(u)=\max\big\{0, \max\{\dfrac{\la x^*_i, u\ra}{b_i}:
i=1,\ldots, q\}\big\}.
\end{eqnarray*}
\end{proposition}
\noindent{\bf Proof. }Let $x^* \in X^*$.  By the definition of $F$ we have that $F$ is a closed convex set.  Since $0\in F$, by Corollary \ref{c1}, one has that $d\in F_\infty$ if and only if
\begin{equation*}
td\in F\mbox{ for all }t>0,
\end{equation*}
or
\begin{equation*}
\la x^*_i, d\ra \leq \dfrac{b_i}{t} \mbox{ for all }i=1,\ldots, q \mbox{ and for all }t>0.
\end{equation*}
This is equivalent to
\begin{equation*}
\la x^*_i, d\ra \leq 0\mbox{ for all }i=1,\ldots, q.
\end{equation*}
The formula for $F_\infty$ has been verified.

Since $b_i>0$ for all $i=1,\ldots, q$, one has $0\in \mbox{int }F$, and
\begin{align*}
\rho_F(u)&=\inf\{t>0: u\in tF\}\\
&=\inf\{t>0: \la x^*_i,u\ra\leq tb_i \mbox{ for all }i=1,\ldots, q\}\\
&=\inf\{t>0: \dfrac{\la x^*_i,u\ra}{b_i}\leq t \mbox{ for all }i=1,\ldots, q\}\\
&=\max\big\{0, \max\{\dfrac{\la x^*_i, u\ra}{b_i}:
i=1,\ldots, q\}\big\}.
\end{align*}
The proof is now complete. $\h$ \vspace*{0.05in}

The following corollary provides a formula for computing the
subdifferential of $\rho_F(u)$ when $u\notin F_\infty$. The formula
allows us to determine a subgradient $x^*_k$ in Theorem
\ref{algorithm} when $F$ is the polyhedron dynamic (\ref{pdy}) in
the case $x_k\notin\Omega_i-F_\infty$.

\begin{corollary} Consider the set $F$ defined in Proposition \ref{poly}. For $u\in X$, define
\begin{equation*}
J(u)=\{i=1,\ldots,p: \dfrac{\la x^*_i,u\ra}{b_i} =\rho_F(u)\}.
\end{equation*}
For any $u\notin F_\infty$, one has $J(u)\neq \emptyset$ and
\begin{equation}\label{submin}
\partial\rho_F(u)=\co \{\dfrac{x^*_i}{b_i}: i\in J(u)\}.
\end{equation}
\end{corollary}
\noindent{\bf Proof. }Since $u\notin F_\infty$, by Proposition \ref{poly}, there exists $j=1,\ldots,q$ such that $\la x^*_j,u\ra >0$. This implies $\rho_F(u)>0$ and $J(u)\neq \emptyset$. Using the formula for computing $\rho_F$ in Proposition \ref{poly} and the well-known formula for computing subdifferential of ``max" function, see e.g. \cite[Corollary~4.3.2]{hl} (which also holds in this setting), one obtains (\ref{submin}). $\h$ \vspace*{0.05in}

The following example shows an application of the subgradient method.

\begin{example}{\rm Let $F=\{(x,y)\in\R^2: x\leq 1, y\leq 1\}$ where $x_1^* = (1,0)$ and $x_2^* = (0,1)$.  Then $F_\infty=\{(x,y):
x\leq 0, y\leq 0\}$. Define the target sets $\Omega_i$ by
two-tuples, that is, $\Omega_i=\{(a_i,b_i)\}$ for $i=1,\ldots,m$ and
let $\Omega= \B(c;r)$. Then $\rho_F((u_1,u_2))=\max\{u_1, u_2, 0\}$,
and
$\mathcal{T}^F_{\Omega_i}((u_1,u_2))=\rho_F((a_i,b_i)-(u_1,u_2))=\max\{0,
a_i-u_1, b_i-u_2\}$. Suppose that we have found $x_k=(x_{1k},
x_{2k})$ from the subgradient algorithm. Then we compute
$\mathcal{T}^F_{\Omega_i}(x_k)=\max\{0, a_i-x_{1k}, b_i-x_{2k}\}$
and $\mathcal{T}(x_k)=\max\{\mathcal{T}^F_{\Omega_i}(x_k):
i=1,\ldots, m\}$. Fix $i\in I(x_k)$. If $x_{1k}\geq a_i$ and
$x_{2k}\geq b_i$, then choose $x^*_k=0$. Otherwise, $x_k\notin
\Omega_i-F_\infty$. If $a_i-x_{1k}\geq b_i-x_{2k}$, then choose
$x^*_k=(-1,0)$. In the case $a_i-x_{1k}< b_i-x_{2k}$, choose
$x^*_k=(0,-1)$. The Euclidean projection to $\Omega$ can be easily
obtained. Thus we can determine $x_{k+1}$, and the algorithm is now
explicit.}
\end{example}

\begin{figure}[h]
\vspace{-1.5cm} \hspace{-1.3cm}
\begin{minipage}{2in}
\vspace{0.3in}
\includegraphics[width=3.0in]{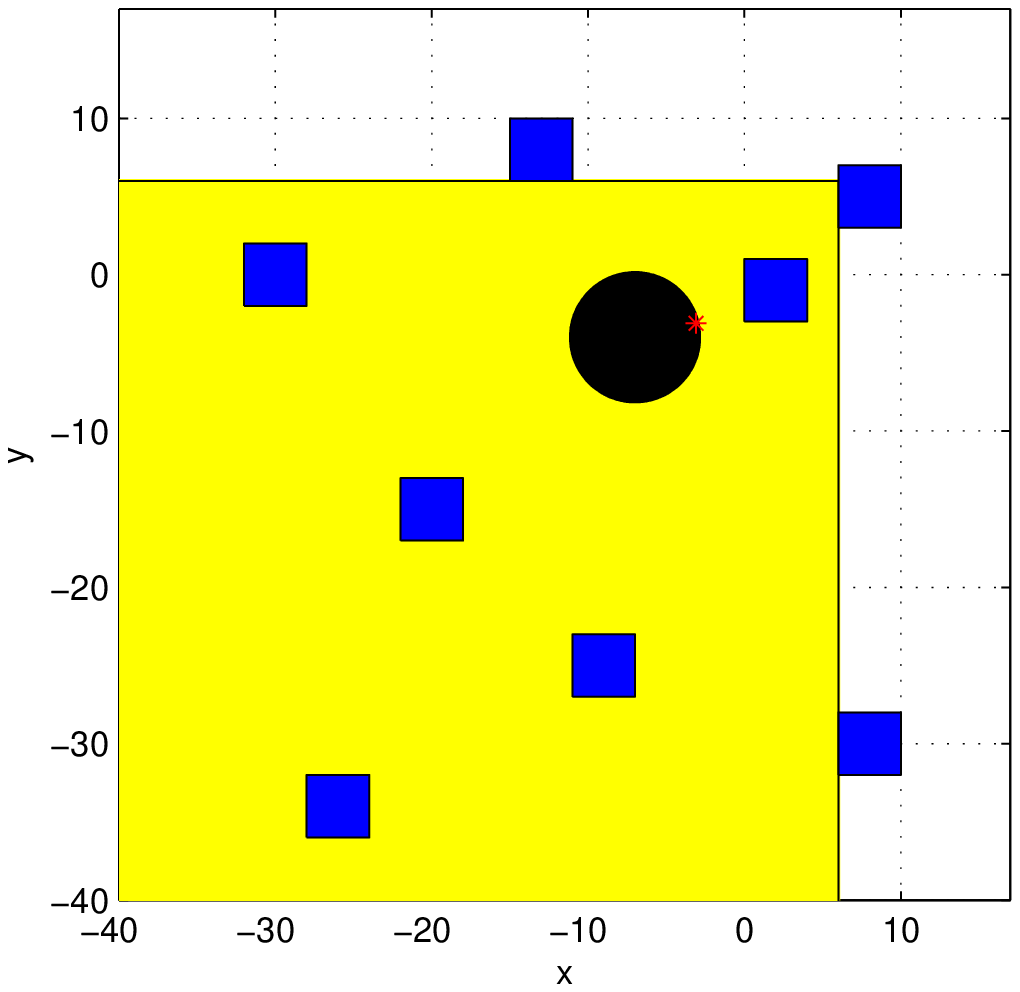}\\
\end{minipage}
~\hfill~ \hspace{1.5cm}
\begin{minipage}[t]{0.54\textwidth}
{\small
\begin{tabular}{|c|c|c|}
\hline
\multicolumn{3}{|c|}{MATLAB RESULT} \\
\hline
$k$ & $x_k$ & $V_k$ \\
\hline
1& (-7.00000, -4.00000)&13.00000\\
10& (-4.17103, -4.00000)& 10.28214\\
100& (-3.10201, -3.10240)& 9.10252\\
1,000&(-3.10210, -3.10200)& 9.10216\\
10,000& (-3.10210, -3.10201)&9.10211\\
100,000& (-3.10209, -3.10208)& 9.10209\\
200,000& (-3.10208, -3.10208)&9.10208\\
300,000& (-3.10208, -3.10208)& 9.10208\\
400,000&(-3.10208, -3.10208)&9.10208\\
500,000& (-3.10208, -3.10208)&9.10208\\
\hline
\end{tabular}
}\end{minipage}
\vspace{-0.5cm}
\caption{A smallest intersecting ball problem to squares with an unbounded dynamic.}\label{Fig1}
\end{figure}
\vspace{-0.5cm}
\begin{example}{\rm Let $F$ be the same set as in the previous example. Let $\Omega_i$ be a square of right
position $[\omega_{1i}-r_i, \omega_{1i}+r_i]\times
[\omega_{2i}-r_i,\omega_{2i}+r_i]$ with center
$c_i=(\omega_{1i},\omega_{2i})$ and short radius $r_i$ for
$i=1,\ldots,m$. Define the constraint $\Omega$ to be a nonempty
closed bounded convex set.  In this case, $$\Omega_i - F_\infty =
\{(x_1, x_2): \ x_1 \geq \omega_{1i} -r_i, x_2 \geq \omega_{2i} -r_i
\}$$ and the projection $\bar\omega$ can be represented by $c_i$ and $r_i$.

The sequence $\{x^*_k\}$ in Theorem \ref{algorithm} can be determined as follows
{\small
\begin{eqnarray*}
x^*_k=\left\{\begin{array}{ll}
(-1, 0),\qquad & \; \mbox{ if } x_{1k} < \omega_{1i} -r_i \mbox{ and }  x_{2k} - \omega_{2i} \geq x_{1k} - \omega_{1i},  \\\\
(0, -1),\qquad & \; \mbox{ if } x_{2k} < \omega_{2i} -r_i \mbox{ and } x_{2k} - \omega_{2i} < x_{1k} - \omega_{1i},  \\\\
(0,0), &\mbox{{otherwise}}
\end{array}\right.
\end{eqnarray*}}
where $i\in I(x_k)$. The sequence $\{V_k\}$ can be determined based on the fact that  $i\in I(x_k)$, so we have
 {\small\begin{eqnarray*}
\mathcal{T}^F_{\Omega_i}((x_{1k},x_{2k}))=\left\{\begin{array}{ll}
(\omega_{1i}-r_i) - x_{1k},\qquad & \; \mbox{ if } x_{1k} < \omega_{1i} -r_i \mbox{ and } x_{2k} - \omega_{2i} \geq x_{1k} - \omega_{1i},  \\\\
(\omega_{2i}-r_i) - x_{2k},\qquad & \; \mbox{ if } x_{2k} < \omega_{2i} -r_i \mbox{ and } x_{2k} - \omega_{2i} < x_{1k} - \omega_{1i},  \\\\
0, &\mbox{{otherwise}}.
\end{array}
\right.
\end{eqnarray*}}

Consider the target sets $\Omega_i, i=1, 2, \ldots, 8$, to be the
right squares with centers $(-30, 0)$, $(-26, -34)$, $(-20, -15)$,
$(-13, 8)$, $(-9, -25)$, $(2, -1)$, $(8, -30)$ and $(8, 5)$ with
radii $r_i=2$ for all $i=1, 2, \ldots, 8$, and  the constraint set
$\Omega$ is the disk with center $(-7, -4)$ and radius 4. Choose
$\alpha_k=1/k$ that satisfies Theorem \ref{algorithm} and the
starting point $(-7, -4)$. An optimal solution and the optimal value
are  $\ox\approx (-3.10208, -3.10208)$ and $\Bar{V}\approx 9.10208$;
see Figure~\ref{Fig1}.}
\end{example}

\section{Concluding Remarks}
\label{s:Concluding}

This paper is a part of our project involving \emph{set facility
location problems}. The main idea is to consider a much broader
situation where singletons in the classical models of facility
location problems are replaced by sets. After developing new results for general and generalized differentiation properties of minimal time functions with unbounded dynamics, we are able to
provide new applications to the smallest intersecting ball problem with convex target sets.

Our next goal is to study the smallest intersecting ball problem with unbounded dynamics and nonconvex target sets. To achieve this goal, we need to develop nonconvex generalized differentiation theory for minimal time functions generated by unbounded dynamics. We foresee the potential of success of
this future work.

\end{document}